\documentclass{amsart}
\usepackage{amssymb}
\usepackage{amsmath}
\usepackage{amstext}
\usepackage{amsfonts}
\usepackage{amsthm}
\usepackage{mathrsfs}
\usepackage{centernot}
\usepackage{enumitem}
\usepackage{graphicx}

\newcommand*{\Scale}
[2][4]{\scalebox{#1}{$#2$}}
\newtheorem{theorem}{Theorem}

\newtheorem*{theorem*}{Lemma}

\newtheorem{remark}{Remark}

\begin{document}

\title[Refinement of the prime geodesic theorem]{On Koyama's refinement of the prime geodesic theorem}
\author{Muharem Avdispahi\'{c}}
\address{University of Sarajevo, Department of Mathematics, Zmaja od Bosne 33-35, 71000 Sarajevo, Bosnia and Herzegovina}
\email{mavdispa@pmf.unsa.ba}
\maketitle

\begin{abstract}
We give a new proof of the best presently known error term in the prime
geodesic theorem for compact hyperbolic surfaces, without the assumption of
excluding a set of finite logarithmic measure. Stronger implications of the
Gallagher-Koyama approach are derived yielding to a further reduction of the
error term outside a set of finite logarithmic measure.
\end{abstract}

\section{Introduction}

Let $\Gamma \subset PSL\left( 2,%
\mathbb{R}
\right) $ be a strictly hyperbolic Fuchsian group acting on the upper
half-plane $\mathcal{H}$ equipped with the hyperbolic metric. The quotient
space $\Gamma \setminus \mathcal{H}$ can be identified with a compact
Riemann surface $\mathcal{F}$ of a genus $g\geq 2$. The object of our
attention is the asymptotic behaviour of the summatory von Mangoldt function
\begin{equation*}
\psi _{\Gamma }\left( x\right) =\sum_{\substack{ P,k  \\ N\left( P\right)
^{k}\leq x}}\log N\left( P\right)
\end{equation*}%
where the sum is taken over primitive hyperbolic conjugacy classes $P$ in $%
\Gamma $ (prime geodesics on $\mathcal{F}$), $N\left( P\right) =\exp \left(
\text{length}\left( P\right) \right) $ is the norm of a class $P$ and $k$
runs through positive integers.

In the recent paper \cite{K}, Shin-ya Koyama
studied the existence of a subset $E$ in $\mathbb{R}_{\geq 2}$ with finite
logarithmic measure such that
\begin{equation*}
\psi _{\Gamma }\left( x\right) =x+\sum_{\frac{3}{4}<\rho <1}\frac{x^{\rho }}{%
\rho }+O\left( x^{\frac{3}{4}}\left( \log \log x\right) ^{\frac{1}{4}%
+\varepsilon }\right) \ \ \left( x\rightarrow \infty,x\notin E\right)\text{.}
\end{equation*}%
Here and in the sequel, $\rho$ denotes zeros of the Selberg zeta function $%
Z_{\Gamma}$. It is known that the complex zeros of $Z_{\Gamma}$ are of the
form $\rho = \frac{1}{2}\pm \mathrm{i}\gamma$ and that $Z_{\Gamma}$ has
finitely many real zeros, all laying in the interval $\left[0,1\right]$.

\noindent Koyama was motivated by Gallagher's \cite{G} approach to the prime
number theorem under Riemann hypothesis.

We give a new proof of the following sharper result (cf. \cite{R}, \cite{B}%
).

\begin{theorem}
\label{T1}
\begin{equation*}
\psi _{\Gamma }\left( x\right) =x+\sum_{\frac{3}{4}< \rho <1}\frac{x^{\rho }%
}{\rho }+O\left( x^{\frac{3}{4}}\right) \ \ \left( x\rightarrow \infty
\right) \text{.}
\end{equation*}
\end{theorem}

We observe that the analogue is also valid for higher dimensional hyperbolic
manifolds with cusps. Applying the Gallagher-Koyama method, we further
reduce the error term outside a set of finite logarithmic measure.

\begin{theorem}
\label{T2} For $\alpha >0$, there exists a set $H$ of finite logarithmic
measure such that
\begin{equation*}
\psi _{\Gamma }\left( x\right) =x+\sum_{\frac{3}{4}-\varepsilon <\rho <1}%
\frac{x^{\rho }}{\rho }+O\left( \frac{x^{\frac{3}{4}}}{(\log x)^{\alpha }}%
\right)  \ \ \left( x\rightarrow \infty ,\ x\notin H\right) \text{,}
\end{equation*}%
where $\varepsilon > 0$ is arbitrarily small.
\end{theorem}

\section{From Hejhal to Randol}

\noindent \textit{Proof of Theorem \ref{T1}}. We shall take the same
starting point as in \cite{K}, i.e. Hejhal's explicit formula with an error
term for the function $\psi _{1,\Gamma }\left( x\right) =\underset{1}{%
\overset{x}{\int }}\psi _{\Gamma }\left( x\right) dx$ (cf. \cite[Theorem
6.16. on p. 110]{H}):
\begin{equation}  \label{one}
\begin{split}
\psi _{1,\Gamma }\left( x\right) =\alpha _{0}x+\beta _{0}x\log x+\alpha
_{1}+\beta _{1}\log x \\
+F\left( \frac{1}{x}\right)+\frac{x^{2}}{2}+\sum_{\substack{ \rho  \\ %
\left\vert \gamma \right\vert <T}}\frac{x^{\rho +1}}{\rho \left( \rho
+1\right) } +O\left( \frac{x^{2}\log x}{T}\right) \ \left( x\rightarrow \infty \right)
\text{.}
\end{split}%
\end{equation}%
Recall that $F\left( x\right) =\left( 2g-2\right) \overset{\infty }{\underset%
{k=2}{\sum }}\frac{2k+1}{k\left( k-1\right) }x^{1-k}$.

The novelty of our approach consists in integrating \eqref{one} at this
point and then temporarily getting rid of Hejhal's error term. Indeed, the
integration of \eqref{one} firstly yields the explicit formula with an error
term for $\psi _{2,\Gamma }\left( x\right) =\underset{1}{\overset{x}{\int }}%
\psi _{1,\Gamma }\left( x\right) dx$. Now, letting $T\rightarrow \infty $ in
the obtained formula, we end up with
\begin{multline*}
\psi _{2,\Gamma }\left( x\right) =\alpha _{0}^{\prime }x^{2}+\beta
_{0}^{\prime }x^{2}\log x+\alpha _{1}^{\prime }x+\beta _{1}x\log x +\frac{x^{3}}{6}+\beta _{2}\\
+\left( 2g-2\right) \overset{\infty }{\underset{k=2}{\sum }}\frac{2k+1}{k\left( k-1\right) }\frac{x^{2-k}}{\left( 2-k\right)} +\underset{\frac{1}{2}< \rho < 1}{\sum }\frac{x^{\rho +2}}{\rho \left( \rho
+1\right) \left( \rho +2\right) }+\underset{\text{Re}\rho =\frac{1}{2}}{\sum }\frac{%
x^{\rho +2}}{\rho \left( \rho +1\right) \left( \rho +2\right) }\text{.}
\end{multline*}

As usually, to derive the asymptotics of $\psi _{\Gamma }\left( x\right) $
from the asymptotics of $\psi _{2,\Gamma }\left( x\right) $, one introduces
the second-difference operators:
\begin{eqnarray*}
\Delta _{2}^{+}f\left( x\right) &=&f\left( x+2h\right) -2f\left( x+h\right)
+f\left( x\right) \text{ and} \\
\Delta _{2}^{-}f\left( x\right) &=&f\left( x-2h\right) -2f\left( x-h\right)
+f\left( x\right) \text{,}
\end{eqnarray*}%
where $h>0$ is to be determined later.

Since $\psi _{\Gamma }$ is a non-decreasing function, we have
\begin{equation*}
\frac{1}{h^{2}}\Delta _{2}^{-}\psi _{2,\Gamma }\left( x\right) \leq \psi
_{\Gamma }\left( x\right) \leq \frac{1}{h^{2}}\Delta _{2}^{+}\psi _{2,\Gamma
}\left( x\right) \text{.}
\end{equation*}%
We apply $\Delta _{2}^{+}$ to all summands in the explicit formula for $\psi
_{2,\Gamma }\left( x\right) $. E.g., $\Delta _{2}^{+}\left( \frac{x^{3}}{6}%
\right) =xh^{2}+h^{3}$, what gives us $\frac{1}{h^{2}}\Delta _{2}^{+}\left(
\frac{x^{3}}{6}\right) =x+h$, etc.

Applying $\frac{1}{h^{2}}\Delta _{2}^{+}$ to the sum $\underset{\frac{1}{2}<
\rho < 1}{\sum }\frac{x^{\rho +2}}{\rho \left( \rho +1\right) \left( \rho
+2\right) }$, we end up with $\underset{\frac{1}{2}< \rho < 1}{\sum }\frac{x^{\rho }}{\rho }+O\left(
h\right)$.

When dealing with the absolutely convergent series $\underset{\text{Re}%
\left( \rho \right) =\frac{1}{2}}{\sum }\frac{x^{\rho +2}}{\rho \left( \rho
+1\right) \left( \rho +2\right) }$, we take into account that
\begin{equation*}
\frac{1}{h^{2}}\Delta _{2}^{+}\frac{x^{\rho +2}}{\rho \left( \rho +1\right)
\left( \rho +2\right) }=O\left( \min \left( \frac{x^{\frac{1}{2}}}{%
\left\vert \rho \right\vert },\frac{x^{\frac{5}{2}}}{h^{2}\left\vert \rho
\right\vert ^{3}}\right) \right)\text{.}
\end{equation*}%
Thus,
\begin{multline*}
\frac{1}{h^{2}}\Delta _{2}^{+}\underset{\text{Re}\left( \rho \right) =\frac{1%
}{2}}{\sum }\frac{x^{\rho +2}}{\rho \left( \rho +1\right) \left( \rho
+2\right) } =O\left( x^{\frac{1}{2}}\underset{{\scriptsize
\begin{array}{c}
\text{Re}\left( \rho \right) =\frac{1}{2} \\
\left\vert \rho \right\vert <M%
\end{array}%
}}{\sum }\frac{1}{\left\vert \rho \right\vert }\right) +O\left( \frac{x^{%
\frac{5}{2}}}{h^{2}}\underset{{\scriptsize
\begin{array}{c}
\text{Re}\left( \rho \right) =\frac{1}{2} \\
\left\vert \rho \right\vert \geq M%
\end{array}%
}}{\sum }\frac{1}{\left\vert \rho \right\vert ^{3}}\right) \\
=O\left( x^{\frac{1}{2}}M\right) +O\left( \frac{x^{\frac{5}{2}}}{h^{2}M}%
\right) \text{ \ for }M>2\text{.}
\end{multline*}%
We are left to optimize the terms $O\left( h\right) $, $O\left( x^{\frac{1}{2%
}}M\right) $, $O\left( \frac{x^{\frac{5}{2}}}{h^{2}M}\right) $. This is
achieved by choosing $h=x^{\frac{3}{4}}$, $M=x^{\frac{1}{4}}$. All other
ingredients are dominated by $O\left( x^{\frac{3}{4}}\right) $.

The same procedure works in case of $\Delta _{2}^{-}\psi _{2,\Gamma }\left(
x\right)$, i.e., for estimating $\psi _{2,\Gamma }\left( x\right)$ from
below.

So,
\begin{equation*}
\psi _{\Gamma }\left( x\right) =x+\sum_{\frac{3}{4}< \rho <1}\frac{x^{\rho }%
}{\rho }+O\left( x^{\frac{3}{4}}\right)\text{.}
\end{equation*}

\begin{remark}
The error term $O\left( x^{\frac{3}{4}}\right) $ in Theorem \ref{T1}. yields
$O\left( x^{\frac{3}{4}}/\log x\right) $ in the prime geodesic theorem.
Concerning the explicit formula for $\psi_{1,\Gamma}$, one can consult \cite%
{AS}, where a better estimate for the logarithmic derivative of the Selberg
zeta function is established.
\end{remark}

\begin{remark}
The full analogue is valid for higher dimensional hyperbolic manifolds with
cusps. Namely, the error term in the prime geodesic theorem in that setting reads $O\left( x^{\frac{3}{2}d_{0}}\left( \log x\right) ^{-1}\right) $, where
$d_{0}=\frac{d-1}{2}$ and $d$ is the dimension of a manifold \cite[Theorem 1]%
{AG}.
\end{remark}

\section{An application of the Gallagher-Koyama method}

\noindent \textit{Proof of Theorem \ref{T2}}. In estimating $\psi _{\Gamma
}\left( x\right) $, we shall use the explicit formula \eqref{one} and the
relation $\frac{1}{h}\Delta _{1}^{-}$ $\psi _{1,\Gamma }\left( x\right) \leq
$ $\psi _{\Gamma }\left( x\right) \leq \frac{1}{h}\Delta _{1}^{+}$ $\psi
_{1,\Gamma }\left( x\right) $, where $0<h<\frac{x}{2}$ is to be determined
later on. Here, $\Delta _{1}^{+}f\left( x\right) =f\left( x+h\right) -f\left( x\right)
$ and $\Delta _{1}^{-}f\left( x\right) =f\left( x\right) -f\left( x-h\right) $.

Let $\beta >4\alpha +1$. According to \eqref{one} and the relation above, we
have
\begin{multline}
\psi _{\Gamma }\left( x\right) \leq
\frac{1}{h}\int_{x}^{x+h}\psi _{\Gamma }\left( t\right) dt \\
=x+\sum_{\frac{1}{2}<\rho <1}\frac{x^{\rho }}{\rho }+O\left( \log x\right)
+O\left( h\right)+O\left( \frac{x^{2}\log x}{hT}\right)
+\frac{1}{h}\left\vert \sum _{\substack{ \text{Re}(\rho )=\frac{1}{2} \\ \left\vert \gamma \right\vert \leq T}}\frac{\left( x+h\right) ^{\rho
+1}-x^{\rho +1}}{\rho \left( \rho +1\right) }\right\vert \text{.}
 \label{two}
\end{multline}%
Now,
\begin{equation*}
\sum_{\substack{ \text{Re}(\rho )=\frac{1}{2}  \\ \left\vert \gamma
\right\vert \leq T}}\frac{\left( x+h\right) ^{\rho +1}-x^{\rho +1}}{\rho
\left( \rho +1\right) }=\sum_{\substack{ \text{Re}(\rho )=\frac{1}{2}  \\ \left\vert \gamma
\right\vert \leq (\log T)^{\beta }}}\frac{\left( x+h\right) ^{\rho
+1}-x^{\rho +1}}{\rho \left( \rho +1\right) }+\sum_{\substack{ \text{Re}(\rho )=\frac{1}{2}  \\ (\log T)^{\beta
}<\left\vert \gamma \right\vert \leq T}}\frac{\left( x+h\right) ^{\rho
+1}-x^{\rho +1}}{\rho \left( \rho +1\right) }\text{.}
\end{equation*}

For the first sum on the right hand side, we have%
\begin{equation*}
\frac{1}{h}\left\vert \sum_{\substack{ \text{Re}\left( \rho \right) =\frac{1%
}{2} \\ \left\vert \gamma \right\vert \leq (\log T)^{\beta }}}\frac{\left(
x+h\right) ^{\rho +1}-x^{\rho +1}}{\rho \left( \rho +1\right) }\right\vert =O\left( x^{\frac{1}{2}}\sum_{\substack{ \text{Re}\left( \rho \right) =\frac{%
1}{2} \\ \left\vert \gamma \right\vert \leq (\log T)^{\beta }}}\frac{1}{%
\left\vert \rho \right\vert }\right) =O\left( x^{\frac{1}{2}}\left( \log
T\right) ^{\beta }\right) \text{.}
\end{equation*}%
The second sum is to be split into
\begin{equation*}
\sum_{\substack{ \text{Re}\left( \rho \right) =\frac{1}{2} \\ (\log
T)^{\beta }<\left\vert \gamma \right\vert \leq T}}\frac{\left( x+h\right)
^{\rho +1}}{\rho \left( \rho +1\right) }-\sum_{\substack{ \text{Re}\left(
\rho \right) =\frac{1}{2} \\ (\log T)^{\beta }<\left\vert \gamma
\right\vert \leq T}}\frac{x^{\rho +1}}{\rho \left( \rho +1\right) } =\sigma _{\beta ,T}\left( x+h\right) -\sigma _{\beta ,T}\left( x\right)
\text{.}
\end{equation*}%
Let
\begin{equation*}
\Scale[0.9]{D_{Y}^{T}=\left\{ x\in \left[ T,eT\right) :\left\vert
\underset{{\scriptsize \begin{array}{c} \text{Re}\left( \rho \right)
=\frac{1}{2} \\ Y<\left\vert \gamma \right\vert \leq T\end{array}}}{\sum
}\frac{x^{\rho +1}}{\rho \left( \rho +1\right) }\right\vert
>\frac{x^{\frac{3}{2}}}{(\log x)^{2\alpha }}\right\} \text{, }Y<T\text{.}}
\end{equation*}%
By Koyama's argument \cite[p. 80]{K},
\begin{equation*}
Y^{-1}\gg \frac{1}{(\log eT)^{4\alpha }}\underset{D_{Y}^{T}}{\int }\frac{dx}{%
x}=\frac{1}{(1+\log T)^{4\alpha }}\mu ^{\times }D_{Y}^{T}\text{.}
\end{equation*}%
Hence,%
\begin{equation*}
\mu ^{\times }D_{Y}^{T}\ll \frac{(1+\log T)^{4\alpha }}{Y}\text{.}
\end{equation*}

For $x \in \left[e^{n}, e^{n+1} \right)$, let $T=e^{n}$. The error term in \eqref{two} becomes $O\left(x \log x / h \right)$. Let $Y$ take values $Y_{1}=\left( \log T\right) ^{\beta }=n^{\beta }$, $%
Y_{2}=\left( n-1\right) ^{\beta }$, $Y_{3}=e^{n-1}$. Denote $%
E_{n}=D_{Y_{1}}^{T}$, $F_{n}=D_{Y_{2}}^{T}$, $G_{n}=D_{Y_{3}}^{T}$ and $%
E=\cup E_{n}$, $F=\cup F_{n}$, $G=\cup G_{n}$, respectively. We have%
\begin{eqnarray*}
\mu ^{\times }E &\ll &\sum_{n=2}^{\infty }\frac{\left( n+1\right) ^{4\alpha }%
}{n^{\beta }}<\infty \text{, since }\beta >4\alpha +1\text{;} \\
\mu ^{\times }F &\ll &\sum_{n=2}^{\infty }\frac{\left( n+1\right) ^{4\alpha }%
}{\left( n-1\right) ^{\beta }}<\infty \text{ for the same reason;} \\
\mu ^{\times }G &\ll &\sum_{n=2}^{\infty }\frac{\left( n+1\right) ^{4\alpha }%
}{e^{n-1}}<\infty \text{.}
\end{eqnarray*}

Put $H=E\cup F\cup G$. Obviously, $\mu ^{\times }H<\infty $. We take $%
x,x+h\in
\mathbb{R}
_{\geq 2}\setminus H$.

For $x\in \left[ e^{n},e^{n+1}\right) \setminus E_{n}$, $T=e^{n}$, we get%
\begin{equation*}
\sigma _{\beta ,T}\left( x\right) =O\left( \frac{x^{\frac{3}{2}}}{(\log
x)^{2\alpha }}\right) \text{.}
\end{equation*}

\noindent Case I. If $x+h\in \left[ e^{n},e^{n+1}\right) \setminus H$, then we
also have%
\begin{equation*}
\sigma _{\beta ,T}\left( x+h\right) =O\left( \frac{\left( x+h\right) ^{\frac{%
3}{2}}}{(\log \left( x+h\right) )^{2\alpha }}\right) =O\left( \frac{x^{\frac{%
3}{2}}}{(\log x)^{2\alpha }}\right) \text{.}
\end{equation*}

\noindent Case II. If $x+h \in \left[ e^{n+1},e^{n+2}\right) \setminus H$, we
shall express the sum $\sigma _{\beta ,T}\left( x+h\right) $ in the form%
\begin{equation*}
\sigma _{\beta ,T}\left( x+h\right) =\underset{\substack{ \text{Re}\left( \rho \right) =\frac{1}{2} \\
n^{\beta}<\left\vert \gamma \right\vert \leq e^{n+1}} }\sum \frac{\left(
x+h\right) ^{\rho +1}}{\rho \left( \rho +1\right) }
-\underset{\substack{ \text{Re}\left( \rho \right) =\frac{1}{2} \\
e^{n}<\left\vert \gamma \right\vert \leq e^{n+1}} }\sum \frac{\left(
x+h\right) ^{\rho +1}}{\rho \left( \rho +1\right) }\text{.}
\end{equation*}

The first sum is $O\left( \frac{x^{\frac{3}{2}}}{(\log
x)^{2\alpha }}\right)$ because $x+h\notin F_{n+1}=D_{Y_2}^{eT}$.
The second sum
\begin{equation*}
\underset{\substack{ \text{Re}\left( \rho \right) =\frac{1}{2} \\ e^{n}<\left\vert \gamma \right\vert \leq e^{n+1}} }{\sum }\frac{\left(
x+h\right) ^{\rho +1}}{\rho \left( \rho +1\right) }=O\left( \frac{x^{\frac{3%
}{2}}}{(\log x)^{2\alpha }}\right) \text{,}
\end{equation*}
since $x+h\notin G_{n+1}=D_{Y_3}^{eT}$.

So, in both cases, the relation \eqref{two} becomes
\begin{equation*}
\psi _{\Gamma }\left( x\right) \leq x+\sum_{\frac{1}{2} < \rho <1}\frac{%
x^{\rho }}{\rho }+O\left( \log x\right) +O\left( h\right)+O\left( \frac{x\log x}{h}\right) +O\left( \frac{x^{\frac{3}{2}}}{h(\log
x)^{2\alpha }}\right) \text{.}
\end{equation*}

The optimal bound is achieved by $h=\frac{x^{\frac{3}{4}}}{(\log
x)^{\alpha }}$. Thus,%
\begin{equation*}
\psi _{\Gamma }\left( x\right) \leq x+\sum_{\frac{1}{2} < \rho <1}\frac{%
x^{\rho }}{\rho }+O\left( \frac{x^{\frac{3}{4}}}{(\log x)^{\alpha }}\right)
\text{.}
\end{equation*}

The opposite inequality is derived from $\psi _{\Gamma }\left( x\right) \geq
\frac{1}{h}\Delta _{1}^{-} \psi _{1,\Gamma }\left( x\right) $ by the same
procedure. If $\varepsilon >0$ is arbitrarily small, then $\sum_{\frac{1}{2} < \rho <\frac{3}{4}-\varepsilon }\frac{%
x^{\rho }}{\rho }$ is obviously dominated by the error term. This completes the proof.

\end{document}